
\input amstex
\documentstyle{amsppt}
\magnification=1200
\NoBlackBoxes
\TagsOnRight
\hcorrection{0.5truein}
\TagsOnRight
\nologo

\hfuzz=9.5pt
\hsize=12.5cm
\vsize=19.2cm
\parindent=0.5cm
\parskip=0pt
\baselineskip=12pt
\topskip=12pt

\widestnumber\key{BDGW}

\def \qed{\hfill$\square $}
\def \-->{\longrightarrow}
\def \kler{K\"ahler}


\def \dbar{\overline{\partial}}



\def\const{{\operatorname{const}}}
\def\det{{\operatorname{det}}}
\def\deg{{\operatorname{deg}}}
\def\dvol{{\operatorname{dvol}}}
\def\rank{{\operatorname{rank}}}
\def\dim{{\operatorname{dim}}}

\def\End{{\operatorname{End}}}
\def\Hom{{\operatorname{Hom}}}

\def\Vol{{\operatorname{Vol}}}


 \def \tri{{(\Cal E_1,\Cal E_2,\Phi)}}

 \def \xp{X\times\Bbb P^1}
 \def \sue{$SU(2)$-equivariant}
 \def \ci{C^{\infty}}

\leftheadtext{}
\rightheadtext{}

\topmatter
\title
A New look at the vortex equations and dimensional reduction
\endtitle
\author Steven~Bradlow \footnote{Supported in part by the NSF
grant DMS-9303545.\hfill \hfill \hfill}, 
James~Glazebrook \footnote{Supported in part by the NSF
grant DMS-9208182.\hfill \hfill \hfill} 
and Franz~Kamber \footnote{Supported in part by the NSF
grants DMS-9208182, DMS-9504084
and by an A. O. Beckman Research Award from the
University of Illinois.\hfill \hfill \hfill}
\endauthor

\address Department of Mathematics,
University of Illinois,
Urbana, IL 61801 ${}^{1,2,3}$
\endaddress

\address Department of Mathematics,
East. Illinois University,
Charleston, IL 61920 ${}^2$
\endaddress
\email bradlow\@uiuc.edu, glazebro\@math.uiuc.edu, kamber\@math.uiuc.edu
\endemail

\subjclass Primary 58C25 ; Secondary 58A30, 53C12, 53C21, 53C55, 83C05
\endsubjclass
\keywords Vortex equation, projective bundles,
 Hermitian--Einstein equation
\endkeywords

\abstract \nofrills {ABSTRACT.}
In order to use the technique of dimensional reduction, it is usually
necessary for there to be a symmetry coming from a group action. In this paper
we consider a situation in which there is no such symmetry, but in which a
type
of dimensional reduction is nevertheless possible. We obtain a relation
between
the Coupled Vortex equations on a closed Kahler manifold, $X$, and the
Hermitian-Einstein equations
on certain $\Bbb P^1$-bundles over $X$. Our results thus generalize the
dimensional reduction results of Garcia-Prada, which apply when the
Hermitian-Einstein equations are on $X\times\Bbb P^1$.
\endabstract
\endtopmatter
\document

\heading
\S 1. Introduction
\endheading

Dimensional reduction is a technique for studying special solutions to partial
differential equations. The technique is applicable when there is a \it
symmetry\rm\ i.e. a group action, in which case it makes sense to look for \it
invariant\rm\ solutions. The term `dimensional reduction'
refers to the fact that the invariant solutions to the original equation can
be interpreted as ordinary solutions to a related set of equations on the
(lower dimensional) orbit space of the group action.

In this paper we will describe a situation in which an effective dimensional
reduction of equations is possible, even though there is no global group
action. The equations in question are the Hermitian-Einstein (HE) equations
for special metrics on holomorphic bundles.  As a result of our dimensional
reduction procedure, we relate special solutions of the HE equations to
solutions of the so-called Coupled Vortex equations.  These latter equations
are also for special metrics on holomorphic bundles, but involve more data
than
HE equations. The extra data is in the form of prescribed holomorphic
sections.

A dimensional reduction from the Hermitian-Einstein to the Vortex equations is
not new. Indeed in [GP], Garcia-Prada described
just such a relation between the equations. In the situation considered by
Garcia-Prada the vortex equations are
defined on bundles over a closed Kahler manifold, say $X$, the
Hermitian-Einstein equations are on bundles over $X\times\Bbb P^1$, and the
reduction is with respect to the standard $SU(2)$-action on $\Bbb P^1$.

The main difference between Garcia-Prada's results and those described in
this paper is that we have replaced $X\times\Bbb P^1$\ by certain \it
non-trivial\rm\ $\Bbb P^1$--bundles over $X$.  One way to view such a
situation
is as one in which there is no group action but in which the (group) orbits
are
still evident. In fact, the group orbits are now the leaves of a foliation
structure.  Instead of reducing to
an orbit space of a group action, we thus reduce to the leaf space of a
foliation.
Such reduction with respect to a foliation has been investigated in the
context
of foliations on Riemannian manifolds in [GKPS]. Interestingly, despite the
suggestiveness of the foliation aspects of our construction, the techniques we
use actually rely more on the holomorphic aspects, and are thus much closer to
those of Garcia-Prada than to those of [GKPS].

Indeed, the results we describe in this paper show that
the Garcia-Prada techniques readily extend to the case where  $X\times\Bbb
P^1$\ is replaced by a \it projectively flat\rm\ $\Bbb P^1$--bundle over $X$.
More precisely, we assume that the $\Bbb P^1$--bundle has a flat $PU(2)$\
structure. This assumption allows us to extend any $PU(2)$-invariant structure
on $\Bbb P^1$\ to the total space, say $M$, of the $\Bbb P^1$--bundle over
$X$.
In particular,  the \kler\ form of the Fubini-Study metric can be extended to
become a global form on $M$, and the canonical bundle on $\Bbb P^1$\ extends
to
become the relative canonical bundle $K_{M/X}$.  By carefully exploiting these
constructions we are able to identify an interesting class of holomorphic
bundles over $M$\ on which the Garcia-Prada techniques can be made applicable.

In Section 2  we give some background material on the Hermitian-Einstein
equations and the Coupled Vortex equations, and in Section 3 we review the
Garcia-Prada dimensional reduction techniques. In section 4 we introduce our
Projective Flatness Condition and discuss its implications. Our main result is
given in section 5, where it appears as Theorem 5.1. Section 6 contains some
computations for the special case that the projective bundle is the
projectivization of a vector bundle. Finally, in section 7 we indulge in some
discussion of our results, with particular emphasis on possible
interpretations
from the point of view of foliation theory.

\bigskip

\noindent \it Acknowledgements.\rm
The authors are pleased to acknowledge useful conversations with, and
informative comments from Tony Pantev, Robert Friedman, D. Toledo and Andrei
Tyurin. The third named author is grateful for the hospitality and support
provided
by the Erwin Schr\"odinger International Institute
for Mathematical Physics in Vienna and the Mathematics Institute of Aarhus
University.

\heading
\S 2. The Equations
\endheading

The Hermitian-Einstein equations determine special metrics on holomorphic
bundles over \kler\ manifolds (cf. [Ko] or [LT] for a complete review).
Indeed, if $(M,\omega_M)$\  is any closed \kler\
manifold, and $\Cal E$\ is a holomorphic bundle over $M$, then the
Hermitian-Einstein equations for a Hermitian metric $H$\ are

$$ i\Lambda_M F_H=\lambda \bold I_E\ .\tag 2.1$$

\noindent Here $E$\ denotes the smooth complex bundle  underlying $\Cal E$,
$F_H$ is the curvature of the metric connection on $\Cal E$ determined by $H$,
$\Lambda_M F_H$ is the
contraction of $F_H$ with the \kler\ form $\omega_M$,
$\bold I_E\in \Omega^0(\End E)$ is the identity and $\lambda$ is a constant
determined by $\omega_M$, the rank of $\Cal E$, and its degree. Recall that
the
degree of a complex bundle over a \kler\ manifold is defined by
 $$
 \deg( E)=\int_M c_1(E)\wedge\omega_M^{m-1}
= \int_M\Lambda_Mc_1(E)\frac{\omega_M^m}{m}\ , \tag 2.2
 $$
 where $m$ is the dimension of  $M$ and $c_1(E)$ is the first Chern class of
$E$. It follows (by the Chern-Weil
formula for $c_1(E)$) that the $\lambda$\ in (2.1) is given by

$$\lambda = \frac{2\pi}{(m-1)!}\frac{1}{\Vol(M)}\frac{\deg(E)}{\rank(E)}\ .$$

We remind the reader that the $\deg(E)$\ is used to define the \it slope\rm,
$\mu(E)$, of the bundle by
$$\mu(E)=\frac{\deg(E)}{\rank(E)}\ .$$
For a holomorphic bundle $\Cal E$, these definitions extend to coherent
subsheaves, and are used to define the algebro-geometric notion of \it
stability\rm. We record, for the sake of completeness, that a bundle $\Cal E$\
is called stable if $\mu(\Cal E')<\mu(\Cal E)$\ for all coherent subsheaves
$\Cal E'\subset\Cal E$.

\bigskip

\noindent\bf Remark 2.1.\rm\ The Hermitian-Einstein equations have many
interesting features, among which are the following:
\roster
\item They are the \kler\ versions of the anti-self-dual equations. What we
mean by this is that, like the ASD equations, the HE equations can be obtained
as the minimizing condition for the Yang-Mills functional on the space of
unitary connections on a complex vector bundle. Furthermore the  ASD and HE
equations are equivalent in the case where they are both defined, i.e. on
unitary bundles with vanishing first Chern class over closed \kler\ surfaces.
\item The existence of solutions is closely related to the  stability of a
holomorphic bundle. This is known as the Hitchin-Kobayashi correspondence.
\endroster

\bigskip

The Coupled Vortex Equations also determine special metrics, but on a \it
holomorphic triple\rm\ over a \kler\ manifold (cf. [GP], [BGP]). If
$(X,\omega_X)$\ is a closed \kler\ manifold of dimension $n$, then by
definition a holomorphic triple on $X$ is a triple
$\tri$ consisting of two holomorphic vector bundles $\Cal E_1$ and $\Cal E_2$
on $X$
together with a  homomorphism $\Phi:\Cal E_2\-->\Cal E_1$, i.e. an element
$\Phi\in H^0(\Hom(\Cal E_2,\Cal E_1))$. The coupled vortex equations on
$\tri$\
are defined to be
$$\align
  i  \Lambda_X F_{H_1}+\Phi\Phi^\ast&= \tau_1 \bold I_{E_1}\\
  i \Lambda_X F_{H_2}-\Phi^\ast\Phi&= \tau_2\bold I_{E_2} ~.
\tag 2.3\endalign$$

In these equations, $\Phi^\ast$\ is the adjoint of $\Phi$\ with respect to the
metrics on $\Cal E_1$ and $\Cal E_2$, and $\tau_1$ and $\tau_2$ are real
parameters. If the
ranks of the bundles are $r_1$\ and $r_2$\ respectively, and we denote their
degrees by $d_1$\ and $d_2$, then the parameters $\tau_1$ and $\tau_2$\
satisfy
the constraint
$$r_1\tau_1+r_2\tau_2=\frac{2\pi}{(n-1)!}\frac{1}{\Vol(X)}
(\deg E_1+\deg E_2)\ .\tag 2.4$$

\bf Remark 2.2.\rm\ The Coupled Vortex equations may be viewed as a natural
generalization of the abelian vortex equations (cf. [JT]). Originally
introduced in the so-called Ginsburg-Landau theory of superconductivity, these
equations were later generalized to describe special Hermitian metrics on
holomorphic bundles with prescribed holomorphic sections. In that setting the
equations play as important a role as do the Hermitian-Einstein equations in
the study of stable holomorphic bundles.
Most recently, the abelian vortex equations have appeared in Seiberg-Witten
theory, where they emerge as the \kler\ version of the Seiberg-Witten
equations (cf. [W]).

\heading
\S 3. Garcia-Prada Dimensional Reduction
\endheading

Suppose now that $M=X\times\Bbb P^1$.  Following Garcia-Prada in [GP], we can
consider the situation in which the product $\xp$
has an $SU(2)$-action in which $SU(2)$\ acts trivially on $X$ and via the
identification with the homogeneous space $SU(2)/U(1)$ on $\Bbb P^1$. We can
thus consider \it SU(2)-equivariant\rm\ bundles on $M$, where a smooth bundle
$V\-->M$\ is called $SU(2)$-equivariant if there is an action of $SU(2)$ on
$V$
covering the action on $\xp$.  In particular, we can consider smooth bundles
over $M=X\times\Bbb P^1$\ of the form
$$V= \pi_1^* E_1\oplus \pi_1^* E_2\otimes \pi_2^* H^{\otimes 2}\ .\tag 3.1$$

Here $\pi_1$ and $\pi_2$ are the projections from $\xp$  to the first and
second factors, $E_i$ is a smooth vector bundle over $X$, and $H$ is the
smooth
line bundle over $\Bbb P^1$ with Chern class 1. The $SU(2)$-action on $V$\ is
trivial on the $\pi_1^* E_i$\ factors and standard on $H^{\otimes 2}$.

Not all the $SU(2)$-equivariant bundles are of the form in (3.1); in general
they can be direct sums where each summand is of the form  $\pi_1^* E \otimes
\pi_2^* H^{\otimes k}$. The special case represented by (3.1) is chosen with a
view towards the dimensional reduction from the Hermitian-Einstein equations
to
the Coupled Vortex equations. For our purposes, the key feature of such
bundles, given in the next Proposition, has to do with the $SU(2)$-equivariant
\it holomorphic\rm\ structures that they admit. Recall that a holomorphic
bundle can be viewed as an underlying smooth bundle with a holomorphic
structure, and that a holomorphic bundle $\Cal V$ is $SU(2)$-equivariant if it
is $SU(2)$-equivariant as a smooth
bundle and in addition the action of $SU(2)$ on $\Cal V$ is holomorphic.

\proclaim{Proposition 3.1 [Prop. 3.9 in [GP]}
There is a one-to-one correspondence between \sue\ holomorphic vector bundles
$\Cal E$ with  underlying
\sue\ $\ci$\ structure  given by (3.1), and holomorphic extensions of
the form

$$
0\-->\pi_1^*\Cal E_1\-->\Cal E\-->\pi_1^*\Cal E_2\otimes\pi_2^*\Cal
O(2)\-->0 ~,\tag 3.2
$$
 where $\Cal E_1$ and $\Cal E_2$ are the bundles in (3.1), equipped with
 holomorphic structures.
Moreover, every such extension is defined by an element
$\Phi\in\Hom(\Cal E_2,\Cal E_1)$.
\endproclaim

The characterization of the extension classes in the above Proposition depends
on the fact that extensions over $X\times\Bbb P^1$  of the form (3.2) are
parametrized by $H^1(X\times\Bbb P^1, \pi_1^*(\Cal E_1\otimes
\Cal E_2^\ast)\otimes \pi_2^*\Cal O(-2))$,  and that by the Kunneth formula
this is isomorphic to $H^0(X, \Cal E_1\otimes \Cal E_2^*)\otimes H^1(\Bbb
P^1,\Cal O(-2))$.  The result thus follows from the fact that $H^1(\Bbb
P^1,\Cal O(-2))\cong\Bbb C$~.

It follows immediately from Proposition 3.1 that there is a
one-to-one correspondence between extensions over
$X\times\Bbb P^1$\ of the type  given in (3.2), and holomorphic triples $\tri$
on $X$. Thus, at the level
of holomorphic objects,

\bigskip

\it a holomorphic triple on $X$\ can be viewed as the dimensional reduction of
an $SU(2)$-equivariant holomorphic bundle on $X\times\Bbb P^1$\rm.

\bigskip

We now turn to the dimensional reduction at the level of equations, i.e. we
examine the relation between the Hermitian-Einstein equations on an extension
as in (3.2) and the Coupled Vortex equations on the corresponding triple.

In order to define the equations, we need to fix \kler\ structures. We let
$\omega_X$\ be the (fixed) \kler\ form on $X$, and let $\omega_{\Bbb P^1}$\
denote the \kler\ form of the Fubini-Study metric on $\Bbb P^1$\ normalized so
that
$$\Vol(\Bbb P^1)=\int_{\Bbb P^1}\omega_{\Bbb P^1}=1\ .$$

\noindent On $X\times\Bbb P^1$\ there is a natural 1-parameter family of
$SU(2)$-equivariant \kler\ metrics, with \kler\ forms
$$\omega_{\sigma}=\pi_1^*\omega_X+\sigma\pi_2^*\omega_{\Bbb P^1}\ .\tag 3.3$$
Using these \kler\ structures, one finds:

\proclaim{Proposition 3.2 (Proposition 3.11 [GP])} Let $T=\tri$ be a
holomorphic triple. Let $\Cal E$ be the \sue\
holomorphic bundle over $\xp$ associated to $T$, that is, let $\Cal E$\ be
given
as an extension
$$
0\-->\pi_1^*\Cal E_1\-->\Cal E\-->\pi_1^*\Cal E_2\otimes\pi_2^*\Cal
O(2)\-->0~.
$$
 Suppose  that $\tau_1$ and $\tau_2$ are related by (2.4) and  let $\sigma$\
be
chosen to satisfy
$$
 \frac{1}{\sigma}=
\frac{(r_1+r_2)\frac{\tau_1 \Vol(X)}{2\pi}-
\frac{(\deg E_1 +\deg  E_2)}{(n-1)!}}{2r_2 \Vol(X)}=
\frac{\tau_1-\tau_2}{4\pi}
\ .\tag 3.4
$$
Fix the \kler\ metric on $X\times\Bbb P^1$\ such that the
\kler\ form is $\omega_{\sigma}$, as given by (3.3).

Then $\Cal E_1$ and $\Cal E_2$ admit metrics satisfying the Coupled Vortex
equations (i.e. equations (2.3)) if and only if $\Cal E$ admits an
$SU(2)$-invariant metric which satisfies the Hermitian-Einstein equation (i.e.
equation (2.1)) with respect to $\omega_{\sigma}$.
\endproclaim

\noindent\bf Remark 3.3.\rm\  Now by the Hitchin-Kobayashi correspondence, the
existence of a solution to the Hermitian-Einstein equation is equivalent (cf.
[UY], [Do], [NS])to the polystability of the holomorphic bundle $\Cal E$.
There is a similar correspondence for holomorphic triples, which relates the
existence of solutions to the coupled vortex equations to
a notion of stability for a holomorphic triple (cf. [BGP]).
As one might expect, Proposition 3.2 thus has an analog which relates \it
stable\rm\ holomorphic extensions on $X\times\Bbb P^1$\ to \it stable\rm\
triples on $X$ (cf. Theorem 4.1 in [BGP]). It should be noted that, with fixed
\kler\ structures on $X$\ and $\Bbb P^1$, the notion of stability for the
holomorphic triple on $X$\ depends explicitly on a parameter $\tau$, while the
notion of stability for the holomorphic bundles on $X\times\Bbb P^1$\ depends
on the choice of \kler\ structure on $X\times\Bbb P^1$. If the \kler\ form on
$X\times\Bbb P^1$\ is taken to be $\omega_{\sigma}$, then the theorem relates
the stability of $\Cal E$\ to the $\tau$-stability of the corresponding triple
on $X$, where $\tau$\ is the same as $\tau_1$\ in Proposition 3.2.  This
analog
to Proposition 3.2 can be be thought of as a dimensional reduction result for
\it stable holomorphic objects ~\rm(rather than for \it equations\rm).

\bigskip

\heading
\S 4. Non-trivial $\Bbb P^1$-bundles over $X$
\endheading

\subheading{\S 4.1 Dimensional reduction of bundles}

It is natural to view $X\times\Bbb P^1$\ as a special case of
a $\Bbb P^1$-fibration over $X$. One is then lead to consider
whether any of the results of the previous section carry over
to this more general situation. Henceforth we thus replace $X\times\Bbb P^1$
by

$$\Bbb P^1\hookrightarrow M\-->X\  ,$$

\noindent a holomorphic $\Bbb P^1$-bundle over X.

In its most general form this
means that $M$\ has a description as a $\Bbb P^1$--bundle associated to
a principal $PGL(2, \Bbb C)$-bundle.  Denoting the principal $PGL(2, \Bbb
C)$--bundle by $\tilde P$, we can thus write

$$M=\tilde P\times_{PGL(2,\Bbb C)}\Bbb P^1\ .$$

\noindent\bf Remark.\rm\ If $X$\ is a Riemann surface, so that $M$\ is a ruled
surface, then the $PGL(2,\Bbb C)$\ structure
always lifts to $GL(2,\Bbb C)$\ (cf. [GH] or [H]) and $M$\ can be described as
$M=\Bbb P(E)$, where $E\-->X$\ is a rank two holomorphic bundle. If
$\dim_{\Bbb
C}X>1$,
then this is not always possible. In general there is an obstruction to such a
lift, with the obstruction located in $H^3(X,\Bbb Z)$.

\bigskip

Our initial goal is to find a class of holomorphic bundles on $M$\
whose structure is determined by \it data on $X$\rm.  In general there is no
longer an $SU(2)$-action on $M$, so it
does not make sense to study $SU(2)$-equivariant bundles over $M$. Instead,
it turns out that the appropriate replacements for the extensions described
in (3.2) are extensions of the form

$$
0\-->\pi^*\Cal E_1\-->\Cal E\-->\pi^*\Cal E_2\otimes\Cal K_{M/X}^*\-->0
\ ,\tag 4.1
$$
where $\pi:M\-->X$\ is the projection, $\Cal E_i$\ are holomorphic bundles
over $X$, and $\Cal K_{M/X}^*$\ is the dual of the relative canonical bundle
on
$M$.


\bigskip
Notice that if $M=X\times\Bbb P^1$, then $\Cal K^*_{M/X}=\pi_2^*\Cal O(2)$\
and we recover the extensions of the form (3.2). Furthermore, an important
part
of Proposition 3.1 remains true for extensions of this sort, namely that
every such extension is defined by an element $\Phi\in\Hom(\Cal E_2,\Cal
E_1)$.
More precisely,

\proclaim{Proposition 4.1} There is a natural isomorphism

$$H^{0,1}(M,\pi^*(\Cal E_1\otimes\Cal E_2^*)\otimes\Cal K_{M/X})
\cong
H^0(X,\Cal E_1\otimes\Cal E_2^*)\ .\tag 4.2$$

\endproclaim

\demo{Proof}

This follows by the Leray spectral sequence, plus the fact that the direct
image
sheaves of the relative canonical bundle satisfy (cf. [H])

$$R^q\pi_*\Cal K_{M/X}\cong 0 \ \text{if}\ q\ne 1\ ,$$
$$R^1\pi_*\Cal K_{M/X}\cong \Cal O_X \ .\tag 4.3$$
\qed\enddemo

What this says is:

\proclaim{Corollary 4.2} Let $M$ be any holomorphic $\Bbb P^1$-bundle over
$X$.
There is a one to one correspondence between holomorphic triples $(\Cal
E_1,\Cal E_2,\Phi)$\ on $X$\ and extensions
$$0\-->\pi^*\Cal E_1\-->\Cal E\-->\pi^*\Cal E_2\otimes\Cal K_{M/X}^*\-->0~,$$
over $M$.
\endproclaim

\subheading{\S 4.2 Projective flatness Condition}  To proceed further with the
dimensional reduction, i.e. to relate the Hermitian-Einstein and Coupled
Vortex
equations, we need to consider the
\kler\ structures on $M$.  In general, it seems to be a difficult matter to
find a \kler\ metric on $M$\ such that the resulting Hermitian-Einstein
equations on an extension (4.1) will `dimensionally reduce' to the Coupled
Vortex equations on the corresponding triple over $X$. We have, however,
discovered a sufficient condition on $M$ in order for this to be possible.
It remains an interesting open question whether our condition is also
necessary.

The condition we impose on $M$ is the

\proclaim{Projective Flatness Condition}
We assume that $\Bbb P^1\hookrightarrow M\-->X$\ is a flat $PU(2)$-bundle.
\endproclaim

The real significance of the projective flatness condition lies in the fact,
that under this assumption, we are able to construct the following:
\roster
\item a class of bundles which correspond to the
$SU(2)$-equivariant extensions in the Garcia-Prada construction,
\item  a suitable family of \kler\ metrics on $M$, generalizing the family
corresponding to the \kler\ forms given in (3.3), and
\item explicit descriptions of the forms representing the extension classes
in (4.1).
\endroster
In short, we can construct bundles which are
are ``locally equivariant'' and retain enough symmetry so that the key
constructions required for dimensional reduction still apply.

Before describing these constructions, it is useful to elaborate somewhat on
the projective flatness condition. There are three equivalent ways of viewing
this condition:
\roster
\item Using the definition of a flat
bundle in terms of locally constant transition functions, the projective
flatness condition means that we can choose a cover $\{U_i\}$\ for $X$,
and local trivializations $M|_{U_i}\cong U_i\times\Bbb P^1$,  such that

$$M= (\coprod_i{U_i\times\Bbb P^1})/g_{ij}\ ,\tag M1$$

\noindent where the transition functions
$$g_{ij}:U_i\cap U_j\-->PU(2)~, $$
\noindent  are locally constant maps. In particular, the coordinate
transformations from $U_j\times\Bbb P^1$\ to
$U_i\times\Bbb P^1$, given by $(x,\lambda)\mapsto(x,g_{ji}\lambda)$, are
holomorphic maps. The description (M1) thus describes $M$\ as a complex
manifold.

\item  Using the cover $\{U_i\}$\ and transition functions $g_{ij}$, we can
construct
a flat principal $PU(2)$--bundle. This is the bundle $\tilde P$\ in the
description given at the beginning of \S 4.1. Hence we see that a second way
to
formulate the projective flatness condition, is to say that

$$M=\tilde P\times_{PU(2)}\Bbb P^1\ ,\tag M2$$
\noindent where $\tilde P$\ is a \it flat\rm\  principal $PU(2)$--bundle.
\item Finally, we can use the fact that a flat bundle $\tilde P$\ can be
associated to a representation

$$\rho:\pi_1(X)\-->PU(2)\ .$$

This gives a description of a projectively flat $M$\ as

$$M=\tilde X\times_{\rho}\Bbb P^1\ ,\tag M3$$

where $\tilde X$\ is the universal cover of $X$, acted on by
$\pi_1(X)$\ via deck transformations.
\endroster

One way that such flat $PU(2)$--bundles arise, is as the projectivization of
rank two holomorphic vector bundles. If $E\-->X$\ is a  vector bundle, then
the
 projective bundle $\Bbb P(E)$\ is defined as follows. Let $E$\ be defined by
local trivializations $E|_{U_i}\cong U_i\times\Bbb C^2$, where  $\{U_i\}$\ is
an open cover of $X$, and  transition functions
$$\hat{g}_{ij}:U_i\cap U_j\--> GL(n+1,\Bbb C)\ ,$$
i.e.
$$E= (\coprod_i{U_i\times\Bbb C^{n+1})/\hat{g}_{ij}}\ .$$

Then  $\Bbb P(E)$\ is defined as

$$\Bbb P(E)= (\coprod_i{U_i\times\Bbb P^n})/g_{ij}\ ,$$
where the transition functions $g_{ij}:U_i\cap U_j\-->PGL(n,\Bbb C)$\ are
obtained from the $\hat{g}_{ij}$\ by the projection from $GL(n+1,\Bbb C)$\ to
$PGL(n,\Bbb C)$. The proof of the following well known fact can be found in
[Ko].

\proclaim{Lemma 4.3}Let $E$\ be a rank $(n+1)$\ holomorphic vector bundle, and
let $H$\ be a Hermitian metric on $E$. If the corresponding metric connection
is projectively flat, i.e. if its curvature $F_H$\ is given by
$$F_H=\alpha\bold I\ ,$$
where $\alpha$\ is a 2-form on $X$\ and $\bold I\in\Omega^0(EndE)$\ is the
identity map, then $\Bbb P(E)$\ has a flat $PU(n+1)$--structure.
\endproclaim

If $X$\ is a Riemann surface, i.e. $\dim_{\Bbb C}X=1$, then
the flat $PU(n)$--bundles which arise in this way are characterized by the
following:

\proclaim{Proposition 4.4}Let $E\-->X$\ be a rank $(n+1)$\ holomorphic vector
bundle over a closed Riemann surface. Then $E$\ admits a projectively flat
Hermitian metric if and only if $E$\ is polystable, i.e. if and only if
$$E=\bigoplus_i{E_i}\ ,$$
where each summand $E_i$\ is a stable holomorphic bundle with
$$\frac{\deg(E_i)}{\rank(E_i)}=\frac{\deg(E)}{\rank(E)}\ .$$

In particular, if $n=1$, then $\Bbb P(E)$\ has a flat $PU(2)$--structure if
and
only if either
\roster
\item $E$\ is a stable rank 2 holomorphic bundle, or
\item $E=L_1\oplus L_2$, where $L_1$\ and $L_2$\ are holomorphic line bundles
with $$\deg(L_1)=\deg(L_2)=\frac{1}{2}\deg (E)\ .$$
\endroster
\endproclaim

\demo{Proof} On a Riemann surface, the condition $F_H=\alpha\bold I$\ is
equivalent (cf. [Ko]) to the Hermitian-Einstein condition $\Lambda
F_H=\lambda\bold I$.
The result thus follows from the Hitchin-Kobayashi correspondence (or, in this
case, the theorem of Narasimhan and Seshadri).\qed
\enddemo

If $\dim_{\Bbb C}X>1$, then the projective flatness condition is in general
stronger than the Hermitian-Einstein condition. However, if $X$\ has a \kler\
metric with \kler\ form $\omega$, and if the Chern classes of $E$\ satisfy the
Bogomolov relation
$$\int_X ((r-1)c_1(E)^2-2rc_2(E))\wedge\omega^{n-2}=0\ ,$$
(where $r=\rank(E)$) then the Hermitian-Einstein condition is again equivalent
to the projective flatness of the full curvature tensor (cf. [LT], also [Ko],
p.114). We thus get:

\proclaim {Proposition 4.5} Let $(X,\omega)$\ be a closed \kler\ manifold, and
let $E\-->X$\ be a rank two holomorphic bundle over $X$. Then $\Bbb P(E)$\ has
a flat $PU(2)$--structure if and only if either
\roster
\item $E$\ is stable (with respect to $\omega$) and
$$\int_X ( c_1(E)^2-4c_2(E))\wedge\omega^{n-2}=0\ ,$$
or
\item $E=L_1\oplus L_2$, where $L_1$\ and $L_2$\ are holomorphic line bundles
with
$$c_1(L_1)=c_1(L_2)\ .$$
\endroster
\endproclaim

\demo{Proof} Part (1) follows immediately from the Bogomolov relation. For
part
(2), we use the fact that the first Chern classes are represented by forms of
type $(1,1)$. They can thus be expressed as $c_1(L_i)=\alpha_i\omega+\beta_i$,
where (for $i=1,2$) $\alpha_i\in \Bbb R$, and the $\beta_i$
are primitive forms of type $(1,1)$, i.e. $\int_X (
\beta_i\wedge\omega^{n-1})=0$. Furthermore, since $E$\ is polystable, it
follows that $\alpha_1=\alpha_2$. The Bogomolov relation thus reduces, in this
case, to
$$\int_X (\beta_1-\beta_2)\wedge(\beta_1-\beta_2)\wedge\omega^{n-2}=0\ .$$

\noindent Using the Hodge-Riemann bilinear relations for real primitive (1,1)
forms (cf. [GH], p.207), we see that the left hand side of this expression is
strictly negative unless $(\beta_1-\beta_2)=0$.\qed
\enddemo

\noindent \bf Remark.\rm\ It should be noted that not all flat $PU(2)$-bundles
arise in this way. As pointed out in the Remark at the beginning of \S 4.1,
there is an obstruction to realizing a $\Bbb P^1$-bundle as $\Bbb P(E)$,
where
$E$\ is a vector bundle.

\bigskip

\subheading{\S 4.3 Equivariant bundles}

Given a $PU(2)$-equivariant holomorphic bundle $\Cal V\-->\Bbb P^1$, we can
use
the description of $M$ as $M=\tilde P\times_{PU(2)}\Bbb P^1$\ to construct
a holomorphic bundle
$$\tilde\Cal V=\tilde P\times_{PU(2)}\Cal V \-->
M=\tilde P\times_{PU(2)}\Bbb P^1\ .\tag 4.4$$
We call $\tilde \Cal V$\ the \it extension of $\Cal V$\ to $M$\rm\ .  In fact
this
construction applies equally well if $PU(2)$\ is replaced by $PGL(2,\Bbb C)$,
and does not require $\tilde P$\ to have a flat structure. If however, $\tilde
P$\ does have a flat structure, so that descriptions (M1-3) apply, then we
get:

\proclaim {Lemma 4.6} Let $\tilde P$ be a flat $PU(2)$-bundle, described as
above by either the representation $\rho:\pi_1(X)\-->PU(2)$\ or the
locally constant transition functions $g_{ij}:U_i\cap U_j\-->PU(2)$.
Let  $\Cal V\-->\Bbb P^1$\ be a $PU(2)$-equivariant holomorphic bundle.

Then the extension of $\Cal V$\ to $M$\ can be described in the following two
equivalent ways:
\roster
\item
$$\tilde\Cal V=\tilde X\times_{\rho}\Cal V\ ,\tag 4.5a$$
\item
$$\tilde\Cal V= (\coprod_i{U_i\times\Cal V})/g_{ij}\ .\tag 4.5b$$
\endroster
\endproclaim

\noindent\bf Remark 4.7.\rm\ Using the above construction, we can extend any
invariant
structure on $\Cal V$\ to $\tilde\Cal V$. In particular,
an invariant Hermitian bundle metric on $\Cal V$\ extends to define a
Hermitian
metric on $\tilde\Cal V$.

\bigskip

When $\Cal V$\ is a line bundle, we can describe the $PU(2)$-equivariant
bundles over $\Bbb P^1$. They are precisely the $SU(2)$-equivariant bundles
on which the element $-I\in SU(2)$\ acts as the identity. The $PU(2)$--action
is then obtained by lifting to an $SU(2)$-action. We conclude therefore that
the $PU(2)$-equivariant holomorphic line bundles on $\Bbb P^1$\ are precisely
the even powers of the tautological bundle $\Cal O_{\Bbb P^1}(-1)$, i.e. are
the bundles $\Cal O_{\Bbb P^1}(2k)$, for any integer $k$. We denote the
corresponding extensions to $M$\ by the notation $\tilde\Cal O(2k)$. For our
purposes, the most important special case is given by $\tilde\Cal O(-2)$.

\proclaim{Lemma 4.8} Let $M$\ be a projectively flat $PU(2)$--bundle, as
above.
Let $K_{M/X}$\ be the relative canonical bundle. Then we can identify
$$K_{M/X}=\tilde\Cal O(-2)\ .$$

\endproclaim


\bigskip

\subheading{\S 4.4 Equivariant metrics}

If $M$\ is projectively flat (and thus has descriptions as in (M1-3)), then
we can construct a 1-parameter family of \kler\ metrics which generalizes the
construction used in [GP]. This is made possible by the fact that $\Bbb P^1$\
has an $SU(2)$-equivariant (and thus a $PU(2)$-equivariant) \kler\ metric,
viz.
the Fubini-Study metric. The corresponding \kler\ form, which we denote by
$\omega_{\Bbb P^1}$, is thus a $PU(2)$-invariant form of holomorphic type
$(1,1)$\ on $\Bbb P^1$. It therefore extends to a form of type $(1,1)$\ on
$M$.
This form, which we denote by $\tilde\omega_{\Bbb P^1}$\ can be described in
two ways (corresponding to the two descriptions (M3) and (M1)):
\roster
\item Let $\pi_2^*\omega_{\Bbb P^1}$\ be the pull-back of $\omega_{\Bbb P^1}$\
to $\tilde X\times \Bbb P^1$, and let
$$q:\tilde X\times \Bbb P^1\-->\tilde X\times_{\rho}\Bbb P^1~,$$
be the quotient map. Then $\tilde\omega_{\Bbb P^1}$\ is the form on $M$
such that $q^*\tilde\omega_{\Bbb P^1}=\pi_2^*\omega_{\Bbb P^1}$.
\item Let  $\omega^{(i)}_{\Bbb P^1}$\ be the pull-back of $\omega_{\Bbb P^1}$\
to $U_i\times \Bbb P^1$, where $\{U_i\}$\ is the cover of $X$\ used in
the description (M1). Because of the $PU(2)$-invariance of $\omega_{\Bbb
P^1}$,
and the fact that the transition functions $g_{ij}$\ are locally constant, the
$\omega^{(i)}_{\Bbb P^1}$\ patch together to define a globally defined form.
The
form they define is $\tilde\omega_{\Bbb P^1}$.
\endroster

\noindent Notice that $\tilde\omega_{\Bbb P^1}$\ is closed and restricts to
$\omega_{\Bbb P^1}$\ on the $\Bbb P^1$\ fibers of $M$.
If $\pi:M\-->X$\ is the projection map, and $\omega_X$\ denotes the \kler\
form
on $X$, then we can define the 1-parameter family

$$\omega_{\sigma}=\pi^*\omega_X + \sigma \tilde\omega_{\Bbb P^1}\ .\tag 4.6$$

This is clearly a family of closed, non-degenerate, positive forms of type
$(1,1)$. In fact $\omega_{\sigma}$\ is the \kler\ form for the Hermitian
metric
on $M$\ which (using description (M1)) pulls back to the (weighted) product
metric on each neighborhood $U_i\times\Bbb P^1$.  This can also be described
as the metric which descends from the weighted product metric on $\tilde
X\times\Bbb P^1$.

In the case where
$M=X\times\Bbb P^1$, this is exactly the family of \kler\ forms described in
[GP]. If $M$\ is not a product, but is projectively flat with a metric of the
above sort, then $M$
is `close enough' to the product case for the following
to be true.

\proclaim{Lemma 4.9} Let $f:X\-->R$\ be a smooth function on $X$. Then

$$\int_M\pi^*(f) \dvol_{\sigma}=\sigma \Vol(\Bbb P^1)\int_X f \ \dvol_X\ ,\tag
4.7a$$
where $\dvol_{\sigma}=\frac{\omega_{\sigma}^{n+1}}{(n+1)!}$\ is the volume
element on $M$, and $\dvol_X=\frac{\omega_X^n}{n!}$\ is the volume element on
$X$.  In particular,
\roster
\item
$$\Vol_{\sigma}(M)=\sigma \Vol(X) \Vol(\Bbb P^1)\ ,\tag 4.7b$$
\item if $V\-->X$\ is a complex bundle on $X$, then
$$\deg_{\sigma}(\pi^*V)=n\sigma~\deg(V)\ ,\tag 4.7c$$
where for any bundle, say $W$, on $M$, $\deg_{\sigma}(W)=\int_M
c_1(W)\wedge\omega_{\sigma}^n$\ is the degree with respect to the \kler\ form
in (4.6).
\endroster
\endproclaim

\demo{Proof} Equation (4.7a) is an immediate consequence of the fact that the
metric is locally a product metric on the open sets homeomorphic to
$U_i\times\Bbb P^1$, and the fact that
$$
\frac{\omega_{\sigma}^{n+1}}{(n+1)!}
=\pi^*(\frac{\omega_X^{n}}{n!})\wedge\sigma
\tilde\omega_{\Bbb P^1}\ .
$$
For (4.7b), take $f=1$. For (4.7c), we use the fact that
$c_1(\pi^*V)=\pi^*c_1(V)$, and also the identity
$$n\alpha\wedge\omega_X^{n-1}=(\Lambda_X\alpha)\omega_X^n\ ,$$
where $\alpha$\ is any complex 2-form on $X$. Using equation (4.6) we then
compute
$$\align
\deg_{\sigma} \pi^* {V}
&=n \sigma \int_M \pi^* (c_1(V) \wedge\omega_X^{n - 1}) \wedge
\tilde\omega_{\Bbb P^1}\ \ \\
&= \sigma \int_M \pi^* (n!\Lambda_Xc_1(V))
\dvol_{\sigma}\ .\\
\endalign$$
Applying (4.7a) with $f=n!\Lambda_Xc_1(V)$\ then gives the result.
\qed\enddemo

\bigskip

\subheading{\S 4.5 Explicit description of extension classes}

In order to adapt the dimensional reduction procedure of Garcia-Prada, we need
to have an explicit description of the extension class in (4.1), i.e. we need
explicit representatives for the classes in
$H^{0,1}(M,\pi^*(\Cal E_1\otimes\Cal E_2^*)
\otimes\Cal K_{M/X})$. Using the Projective Flatness condition, we now show
that these can be taken to be of the form $\pi^*\Phi\otimes\tilde\eta$, where
$\Phi\in\Omega^0(X,\Cal E_1\otimes\Cal E^*_2)$\ and
$\tilde\eta$\ is the extension to $\Cal K_{M/X}$\ of a
uniquely determined invariant element $\eta\in \Omega^{(0,1)}(\Bbb P^1,\Cal
O(-2))$.

We begin with a description of $\eta$. Up to a scale factor, this is uniquely
determined by the requirements that
\roster
\item $\eta$\ is $SU(2)$- (and thus also $PU(2)$-) invariant,
\item $\eta$\ represents a generator of $H^{(0,1)}(\Bbb P^1,\Cal
O(-2))\cong\Bbb C$.
\endroster
We can give an explicit description of $\eta$\ if we
use the identification of $\Cal O(-2)$\ with $K_{\Bbb P^1}$\ i.e. with the
holomorphic cotangent bundle.  Let $z$\ be a local coordinate on $\Bbb P^1$,
and use $\{dz\}$\ as a local frame for $K_{\Bbb P^1}$. With respect to this
frame, and with respect to the local coordinate $z$, we get

$$\eta(z)= \eta_0 d\overline{z}\otimes \frac{dz}{(1+|z|^2)^2}\ ,\tag 4.8$$
\noindent where $\eta_0$\ is the undetermined scale factor.
Now if we interpret the Fubini-Study metric as a bundle metric on the
holomorphic tangent bundle $T^{1,0}\Bbb P^1$,
then with respect to the local coordinate $z$, and the local frame ${\partial
\over \partial z}$, this metric is
given by
$$k({\partial \over \partial z},{\partial \over \partial z})=
\frac{1}{(1+|z|^2)^2}\ .\tag 4.9$$
Thus taking duals with respect to this metric, we can
write
$$\eta(z)= \eta_0 d\overline{z}\otimes
( {\partial \over \partial z})^*\ .$$
Or, if $\eta ^*$\ denotes the `conjugate adjoint', i.e. the section of
$\Omega^{(1,0)}(\Bbb P^1,\Cal O(2))$\ obtained from $\eta$ by complex
conjugation on the form part of $\eta$\ and taking the adjoint (with respect
to
the metric induced by $k$) of the section of $\Cal O(-2)$, then
$\eta ^*=\overline{\eta_0} dz\otimes
 {\partial \over \partial z}$. We thus get that

$$\eta\wedge\eta^*= i|\eta_0|^2\omega_{\Bbb P^1}\ .\tag 4.10$$

\noindent With a view towards the next section, we henceforth fix $\eta_0$\
such that
$$|\eta_0|^2=\frac{1}{\sigma}\ ,\tag 4.11$$
\noindent where $\sigma$\ is the weighting factor in the
\kler\ metric on $M$, i.e. in (4.6).

\proclaim{Definition 4.10} We define
$$\tilde\eta\in \Omega^{(0,1)}(M,\tilde\Cal O(-2))\cong
\Omega^{(0,1)}(M,\tilde\Cal K_{M/X}) ~, $$
to be the extension of $\eta$\ (defined by (4.8) and (4.11)).
\endproclaim

\noindent A local calculation confirms that
this is a $\dbar$-closed form.  Furthermore, the cohomology class which it
represents in $H^{(0,1)}(\Bbb P^1,\Cal O(-2))\cong\Bbb C$\ is non-trivial
since, for example, it restricts to
the generator of $H^{(0,1)}(\Bbb P^1,\Cal O(-2))$\ on each $\Bbb P^1$ fiber of
$M$.  We thus get

\proclaim {Proposition 4.11} The isomorphism

$$ H^0(X,\Cal E_1\otimes\Cal E_2^*)\cong
H^{0,1}(M,\pi^*(\Cal E_1\otimes\Cal E_2^*)\otimes\Cal K_{M/X})~, $$

is realized by the map
$$\Phi\-->\pi^*\Phi\otimes\tilde\eta\in
\Omega^{0,1}(M,\pi^*(\Cal E_1\otimes\Cal E_2^*)\otimes\Cal K_{M/X})~.\tag
4.12$$
\endproclaim

\demo{Proof} We need to check that $\pi^*\Phi\otimes\tilde\eta$\ is
$\dbar$-closed (so that the map is well defined at the level of the cohomology
groups), and that the map is injective
(and thus an isomorphism). The first issue is clear, since $\Phi$\ is assumed
to be a holomorphic section of $\Cal E_1\otimes\Cal E_2^*$ and $\tilde\eta$\
is
also $\dbar$-closed.  The injectivity of the map follows from the above
comment
about the non-triviality of the cohomology class of $\tilde\eta$.\qed
\enddemo

The following property of $\tilde\eta$\ will be used in the next section:

\proclaim{Lemma 4.12}Let $\tilde k$\ be the extension to $\tilde\Cal O(2)\cong
 \Cal K^*_{M/X}$\ of the metric $k$\ on $\Cal O(2)$.  Using this metric, and
the
metric induced by it on  $\tilde\Cal O(-2)\cong \Cal K_{M/X}$, we get

$$\tilde\eta\wedge\tilde\eta^*= \frac{i}{\sigma}\tilde\omega_{\Bbb P^1}\ .\tag
4.13$$
\endproclaim

\demo{Proof}  This is a pointwise property. But in local coordinates on a
neighborhood biholomorphic to $U_i\times\Bbb P^1$, this corresponds precisely
to the fact that $\eta\wedge\eta^*= \frac{i}{\sigma}\omega_{\Bbb P^1}$.
\qed\enddemo

\bigskip

\heading
\S 5. Dimensional Reduction of the Hermitian-Einstein equations
\endheading

We are now ready to examine the Hermitian-Einstein equations
for an extension of the form in (4.1) over a projectively
flat PU(2)--bundle $M$, with respect to a \kler\ metric as above.  The setup
we
consider is the following:

Let $M$\ be a projectively flat $PU(2)$--bundle over $X$, with descriptions as
in (M1), (M2), and (M3). Fix a \kler\ metric on $M$\ with \kler\ form
$\omega_{\sigma}$\ as in equation (4.6). Let $\Vol_{\sigma}(M)$\ denote the
volume with respect to this metric (cf. Lemma 4.9), and for any vector bundle
$V\-->M$, set
$$\deg_{\sigma}(V)=\int_M c_1(V)\wedge\omega_{\sigma}^n\ .$$

Consider a holomorphic extension over $M$\ of the form
$$
0\-->\pi^*\Cal E_1\-->\Cal E\-->\pi^*\Cal E_2\otimes\Cal K_{M/X}^*\-->0~,
$$
as in (4.1). The extension class is represented by a holomorphic section in
$$\Omega^{(0,1)}(M, \Hom(\pi^*\Cal E_2\otimes
\Cal K_{M/X}^*,\pi^*\Cal E_1))
\cong\Omega^{(0,1)}(M,\pi^*(\Cal E_1\otimes\Cal E_2^*)\otimes\Cal K_{M/X})\
.$$
Using the isomorphism given in Proposition 4.11, we can take this to be
$$\beta=\pi^*\Phi\otimes\tilde\eta\ .
\tag 5.1$$


\proclaim{Theorem 5.1}  Let $\Cal E$\ be a holomorphic  extension as above,
with extension class $\beta$\ as in (5.1).  Let $\bold h$ be
a Hermitian metric on ${\Cal E}$ given by

$$\bold h=\pi^* h_1 \oplus\pi^* h_2\otimes \tilde k\ ,$$
\noindent where, for $i =1,2$, $h_i$ is a Hermitian metric on ${\Cal E}_i$
and $\tilde k$\ is the Hermitian metric on $\Cal K_{M/X}^*=\tilde \Cal O(2)$\
described in Lemma 4.12.

Fix $\sigma >0$, and set
$$
\lambda =  \frac{2\pi }{n!}\frac{1}{\Vol_{\sigma}(M)}
{\deg_{\sigma}({\Cal E}) \over{\rank ~{\Cal E}}}\ .\tag 5.2
$$

Let parameters $\tau_1$ and $\tau_2$ be given by
$$ \align
\tau_1&= \lambda\\
\tau_2&=
\lambda -
\frac{2\pi}{n!\Vol_{\sigma}(M)} \deg_{\sigma}(\Cal K_{M/X}^*)\ .\tag 5.3
\endalign$$

\medbreak
Then the following are equivalent:
\roster
\item
The metric $\bold h$ satisfies the Hermitian--Einstein equation
$$
i\Lambda_\sigma F_{\bold h} = \lambda  \bold I_{{\Cal E}}~,
$$
where $\Lambda_\sigma$ denotes contraction
against the
K\"ahler form
$$
\omega_{\sigma} =
\pi^*\omega_X+ \sigma\tilde\omega_{\Bbb P^1}\ .
$$

\item
The metrics $h_1$ and $h_2$ satisfy
the coupled vortex equations:
$$
\align
i\Lambda_X F_{h_1} +\Phi  \Phi^{*} &=
\tau_1 \bold I_{{\Cal E}_1}    \\
i\Lambda_X F_{h_2} -\Phi^{*}  \Phi  &=
\tau_2 \bold I_{{\Cal E}_2}\ ,
\endalign
$$
where the adjoint in $\Phi^{*}$ is with respect to the
metrics $h_1$ and $h_2$, and  $\Lambda_X$ denotes contraction
against the K\"ahler form $\omega_X$\ on $X$.

\endroster
\endproclaim

\demo{Proof}
 Because of the projectively flat structure on $M$\ and our choice of
\kler\ structure, the geometry of our situation is \it locally\rm\
indistinguishable from the case of $M=X\times\Bbb P^1$. The proof of this
theorem is thus essentially the same as that of the corresponding
result in [GP].  We  proceed as follows.

Given a Hermitian metric $\bold h$\ on $\Cal E$, we can analyse the
Hermitian-Einstein
tensor $i\Lambda_{\sigma} F_{\bold h}$. With respect to a smooth
orthogonal splitting $\Cal E=\pi^*\Cal E_1\oplus\pi^*\Cal E_2\otimes\Cal
K_{M/X}^*$, we can write
$$
i\Lambda_{\sigma}F_{\bold h} =
i\Lambda_{\sigma}\pmatrix
F_{\bold h_1} -
\beta \wedge \beta^* & D'\beta \\
-D''\beta^* & F_{\bold h_2} - \beta^* \wedge \beta
\endpmatrix ~. \tag 5.4$$
Here $\bold h_1$\ and $\bold h_2$\ are the metrics induced on the sub- and
quotient bundles,
the $F_{\bold h_i}$\ are the curvatures of the corresponding metric
connections,   and

$$D=D'+D'': \Omega^1(\pi^*(\Cal E_1\otimes\Cal E_2^*)\otimes\Cal K_{M/X})
\rightarrow \Omega^2(\pi^*(\Cal E_1\otimes\Cal E_2^*)\otimes\Cal K_{M/X})$$

\noindent is constructed in the standard way from the metrics on
$\pi^*(\Cal E_1)$\ and
$\pi^*(\Cal E_2^*)\otimes\Cal K_{M/X}$. The summands $D'$ and $D''$\ are the
$(1,0)$\
and $(0,1)$\ parts respectively. The ``$^*$'' in $\beta^*$\ denotes the
adjoint
on sections (with respect to the bundle metrics) and conjugation on forms, so
that
$\beta^*\in \Omega^{(1,0)}(M,\pi^*(\Cal E_2\otimes\Cal E_1^*)\otimes\Cal
K^*_{M/X})$.  We can write
$$ \align
F_{\bold h_1}&=\pi^* F_{h_1}\\
F_{\bold h_2}&=\pi^* F_{h_2}+ F_{\tilde k}\otimes\bold I_2\ .\tag 5.5
\endalign$$

\noindent Now the $SU(2)$-invariant metric $k$\ on $\Cal O_{\Bbb P^1}(2)$\
(described in (4.9)) is Hermitian-Einstein with respect to the Fubini-Study
metric on $\Bbb P^1$. Because of our choice of \kler\ structure on $M$\, it
thus follows that $\tilde k$\ is a Hermitian-Einstein metric
on $\Cal K_{M/X}^*$, so that we have
$$i\Lambda_{\sigma} F_{\tilde k}= c=
\frac{2\pi}{n!\Vol_{\sigma}(M)}\deg_{\sigma}(\tilde\Cal O(2)) \ .\tag 5.6$$

Also, with $\beta$\ as in (5.1), we get
$$i\Lambda_{\sigma}(\beta\wedge\beta^*)=
\pi^*\Phi\Phi^*\otimes i\Lambda_{\sigma}(\tilde{\eta}\wedge\tilde {\eta}^*)
\ ,\tag 5.7$$
where the adjoint in $\Phi^*$\ is with respect to the
metrics on $\Cal E_1$\ and $\Cal E_2$, and $\tilde {\eta}^*$\ is determined by
the metric $\tilde k$.
The proof of the theorem thus depends on the following lemma:

\proclaim{Lemma 5.2} With $\tilde\eta,\ \beta$, and $\Lambda_{\sigma}$
as above, we have
\roster
\item $i\Lambda_{\sigma}(\tilde\eta\wedge\tilde\eta^*)=-1$,
\item $\Lambda_{\sigma}D'\beta=0$, and
\item $\Lambda_{\sigma}D''\beta^*=0$.
\endroster
Furthermore, for $i=1,2$, we get

(4)  $\Lambda_{\sigma}\pi^*F_{h_i}= \pi^*\Lambda_X F_{h_i}$.

\endproclaim
\demo{Proof of Lemma 5.2}

Part (1) follows immediately from Lemma 4.12.

For part (2), write

$$D'\beta=D'(\pi^*\Phi)\otimes\tilde\eta + \pi^*\Phi\otimes D'\tilde\eta\ .$$

We are abusing notation here, since the covariant derivatives denoted by $D'$\
are not the same in the three terms in this expression. We can be more precise
if we restrict to local coordinates on a neighborhood biholomorphic to
$U_i\times\Bbb P^1$.  The formula then becomes

$$D'\beta=\pi_1^*D'_{\Cal E_1\otimes\Cal E_2^*}\Phi\otimes\pi_2^*\eta +
\pi_1^*\Phi\otimes \pi_2^*D'_{\Cal O(-2)}\eta\ ,\tag 5.8$$
where the $\pi_i$\ are the projections onto the first and second factors of
$U_i\times\Bbb P^1$. The covariant derivatives are now those corresponding to
the metric connections on $\Cal E_1\otimes\Cal E_2^*$\ and on $\Cal O_{\Bbb
P^1}(-2)$. Furthermore, in these local coordinates,

$$\omega_{\sigma}=\pi_1^*\omega_X+\sigma\pi^*_2\omega_{\Bbb P^1}\ .$$

The rest of the proof is essentially the same as in [GP], with the key points
being:

\noindent (i) The term
$\Lambda_{\sigma}\pi_1^*D'_{\Cal E_1\otimes\Cal E_2^*}\Phi\otimes\pi_2^*\eta$
\ vanishes because the \kler\ form has no contribution from `mixed' terms,
i.e.
terms in $\Omega^{(1,0)}(U_i)\wedge\Omega^{(0,1)}(\Bbb P^1)$, and

\noindent (ii) the term
$\pi_1^*\Phi\otimes \Lambda_{\sigma}\pi_2^*D'_{\Cal O(2)}\eta$
vanishes because
$$\Lambda_{\sigma}\pi_2^*D'_{\Cal O(2)}\eta=
\pi_2^*(\sigma\Lambda D'_{\Cal O(-2)}\eta)=0\ .$$

The proof of part (3) is similar to (2) (alternatively, one can simply observe
that $\Lambda_{\sigma}D''\beta^*=-(\Lambda_{\sigma}D'\beta)^*$).

For part(4), we note that with $\omega_{\sigma} =
\pi^*\omega_X+ \sigma\tilde\omega_{\Bbb P^1}$, we get

$$\Lambda_{\sigma}\pi^*F_{h_i}=\pi^*\Lambda_X F_{h_i}+\sigma
(\pi^*F_{h_i},\tilde\omega_{\Bbb P^1})\ ,\tag 5.9$$

\noindent where the inner product in the second term on the right is on
$\Omega^2(M,\Bbb C)$. A computation in local coordinates verifies that for any
$\alpha\in\Omega^2(X,\Bbb C)$, we get $(\pi^*\alpha,\tilde\omega_{\Bbb
P^1})=0$.\qed

\enddemo

\noindent\bf Remark 5.3 \rm\  The proof of Lemma 5.2 uses the fact that $M$\
is
projectively flat and the special properties of the \kler\ metrics
that we fix on $M$.  In a more general situation, it is not so clear which
(if any)  \kler\ structures are suitable for the carrying out of a
dimensional reduction of the
Hermitian-Einstein equations. It is tempting to speculate that
the properties listed in the lemma should serve as the \it definition\rm\  of
suitable  \kler\ structures.
\bigskip
\newpage
\noindent\it Completion of Proof of Theorem 5.1\rm\

Using (5.2), (5.6), and Lemma 5.2, we thus find that the Hermitian-Einstein
condition on $\bold h$\ is equivalent to

$$
\pmatrix
\pi^*(i\Lambda_X F_{h_1} +\Phi\Phi^*) & 0 \\
0 & \pi^*(i\Lambda_X F_{h_2} - \Phi^*\Phi) +c\otimes\bold I_2
\endpmatrix=\lambda
\pmatrix
\bold I_1 & 0 \\
0 &  \bold I_2
\endpmatrix
\tag 5.10$$

\noindent where $c$\ is as in (5.6).  The results follows from this. \qed
\enddemo

\noindent\bf Remark 5.4 \rm\  By using (5.3), (5.4), and (4.7c) we can compute
$$\align
r_1\tau_1+r_2\tau_2 &= \frac{2\pi}{n! \Vol_{\sigma}(M)}
(\deg_{\sigma}(\Cal E)-\deg_{\sigma}(\Cal K_{M/X}^*))\\
&= \frac{2\pi}{n! \Vol_{\sigma}(M)}
(\deg_{\sigma}(\pi^*(\Cal E_1\oplus\Cal E_2))\\
&= \frac{2\pi}{(n-1)! \Vol(X)}\deg (\Cal E_1\oplus\Cal E_2)\ ,
\endalign$$
thus verifying that the relation given by (2.4) holds. Also, by (5.4) and
(4.7b), we have

$$\tau_1-\tau_2 = \frac{2\pi}{\sigma}
\frac{\deg_{\sigma}(\Cal K_{M/X}^*)}{n! \Vol(X)}\ .\tag 5.11$$
We can compute $\deg_{\sigma}(\Cal K_{M/X}^*)$, using the fact that $\Cal
K_{M/X}^*=\tilde{\Cal O}(2)$:

\proclaim{Lemma 5.5}Let $\tilde{\Cal O}(2k)$\ be the extension to $M$\ of the
bundle $\Cal O_{\Bbb P^1}(2k)$. Then the first Chern class $c_1({\tilde{\Cal
O}(2k)})$\ can be represented by the $(1,1)$-form $2k\tilde\omega_{\Bbb P^1}$,
where $\tilde\omega_{\Bbb P^1}$\ is the extension of the Fubini Study \kler\
form
(as
in \S 4.4). Hence (assuming that the Fubini-Study metric is normalized so that
$Vol(\Bbb P^1)=1$) we get
$$\deg_{\sigma}(\tilde{\Cal O}(2k))=2kn!Vol(X)\ .$$
\endproclaim

\demo{Proof}  The first Chern class of $\Cal O_{\Bbb P^1}(2k)$\ can be
represented by the (1,1)-form $2k\omega_{\Bbb P^1}$. This in turn can be
obtained as $\frac{i}{2\pi}\partial\dbar log(h)$, where $h$\ is an
($SU(2)$-invariant) hermitian metric on $\Cal O_{\Bbb P^1}(2k)$. Since $h$\ is
$PU(2)$-invariant, it extends to a metric $\tilde h$\ on $\tilde{\Cal O}(2k)$.
One can compute that
$$\frac{i}{2\pi}\partial_M\dbar_M log(\tilde h)
=2k\tilde \omega_{\Bbb P^1}\ ,$$
where $\partial_M$\ and $\dbar_M$\ are the holomorphic and anti-holomorphic
parts of the exterior derivative on $M$. It follows immediately from this that
$c_1({\tilde{\Cal O}(2k)})$\ can be represented by $2k\tilde\omega_{\Bbb
P^1}$.
The formula for $\deg_{\sigma}(\tilde{\Cal O}(2k))$\ then follows from the
form
of $\omega_{\sigma}$\ (cf. (4.6)), and by the computation for
$Vol_{\sigma}(M)$\ (cf. (4.7b)).\qed
\enddemo

Applying Lemma 5.5. to $\Cal K_{M/X}^*=\tilde{\Cal O}(2)$, we get that
$\deg_{\sigma}(\Cal K_{M/X}^*)= 2n! \Vol(X)\ $.  Equation (5.11), which
applies
in the case where $M$\ is a flat $PU(2)$-bundle, is therefore the same as the
relation given in Proposition 3.1, which refers to the case where
$M=X\times\Bbb P^1$.

\heading
\S 6.  Special case where $M=\Bbb P(E)$
\endheading

As an interesting special case, we consider the situation in which $M$\
comes from a holomorphic vector bundle with a projectively flat Hermitian
structure. In this case, the parameter computations of the previous section
can
then be carried out quite explicitly. We thus assume that:
\roster
\item $M=\Bbb P(E)$, where $E\-->X$\ is a rank 2 holomorphic bundle, and
\item there is a Hermitian metric $h$\ on $E$\ such that
$$i\Lambda_X F_h = \const.\ \bold I= \frac {\pi\deg(E)}{(n-1)!\Vol(X)}\
\bold I\
.$$
\endroster

Using the description $M=\Bbb P(E)$, we get a canonically defined line bundle
on $M$, namely the tautological line bundle on $\Bbb P(E)$. This bundle, which
we denote by $\Cal O_M(-1)$, is a subbundle of $\Bbb P(E)\times E$\ and
restricts to the tautological line bundle on each $\Bbb P^1$--fiber of $\Bbb
P(E)$.  The main result of this section is :

\proclaim{Proposition 6.1}
$$\align
\deg_{\sigma}({\Cal O}_M(-1))
&= \frac{n \sigma \deg (E)}{2}  -  n! \Vol (X)\tag 6.1a\\
&= \frac{\deg_{\sigma}(\pi^*(\det E))}{2}-  n! \Vol (X)\ .\tag 6.1b
\endalign$$
\endproclaim

\demo{Proof} The proof is a computation, using Equation (2.2) and the
Chern-Weil formula for $c_1(E)$\ in terms
of curvature.  Since ${\Cal O}_M(-1)$ is a subbundle of $\Bbb P(E) \times
\pi^*
E$, the metric $h$ on $E$ induces a fiber metric, say $H$, on ${\Cal
O}_M(-1)$.
Let $F_H$ be the curvature of its corresponding metric connection.
We thus need to compute ${{i}\over {2\pi}}\Lambda_{\sigma}F_H$ at a point of
 $\Bbb P(E)$.

Any point in $\Bbb P(E)$\ can be represented by a unit vector $\xi_0$ in a
fiber $E_{x_0}$ of $E$. As in [Ko] (p. 89)],  we choose a normal  holomorphic
frame $s=(s_1,s_2)$ for $E$ with $s_2(x_0) = \xi_0$. This gives local
coordinates
$(z^1, \ldots, z^n, \xi^1 ,\xi^2)$\ on $E$, where
$(z^1, \ldots, z^n)$ are local coordinates near $x_0$ on $X$, and $(\xi^1
,\xi^2)$\ are the coordinates with respect to the frame $s$\ of a point in the
fiber. If we denote the points in $E$\ by $\xi(z)$, and the corresponding
points in $\Bbb P(E)$\ by $[\xi(z)]$. Then
$$
[\xi(z)] \mapsto (z^1, \ldots, z^n , [\xi^1,\xi^2])~,
$$
gives a local trivialization
$\Bbb P(E) \vert_{U} \cong U \times \Bbb P^1$, and thus defines a local
holomorphic frame for  ${\Cal O}_M(-1)$ .
Using these local coordinates, the computation in [Ko2 (p.90 and I.4)] shows
that
$$
\align
(i F_H)_{[\xi_0]}  &= i \sum R_{i \bar k \alpha \bar \beta}~
\xi_0^i \xi_0^{\bar k}~  dz^{\alpha} \wedge dz^{\bar \beta}
{}~-~(\omega_{\Bbb P^1})_{[0,1]}   \\
&= i h (R(\xi_0), \xi_0)~-~(\omega_{\Bbb P^1})_{[0,1]} ~,
\tag 6.2
\endalign
$$
where $R$ is the curvature of the metric connection for $h$ on $E$.
But , since $h$\ satisfies the Hermitian--Einstein condition, we get that

$$
{i n }~h(R(\xi_0), \xi_0) \wedge \omega^{n-1}_X
=  \frac{\pi \deg (E)}{\Vol(X)(n-1)!}{\vert \xi_0 \vert}^2 ~ \omega_X^n~.
\tag 6.3
$$
Also, using the local coordinate system
on $\Bbb P(E)$ near $[\xi_0]$,
the \kler\ form on $M=\Bbb P(E)$\ is given by $
\omega_\sigma =  \omega_X + \sigma \omega_{\Bbb P^1}$.
Thus, using ${\vert \xi_0 \vert}^2 = 1$, we get
$$
\align
{{(i F_H)} \over {2 \pi}} \wedge \omega_{\sigma}^n
&= {{(i n \sigma)} \over {2 \pi}} h(R(\xi_0),\xi_0)
\wedge \omega_X^{n - 1} \wedge \omega_{\Bbb P^1} ~-~
\omega_X^{n} \wedge \omega_{\Bbb P^1}   \\
&= (\frac{n \sigma\deg (E)}{2 \Vol(X)} - n!)(\frac{\omega_X^{n}}{n!} \wedge
\omega_{\Bbb P^1})\\
&= (\frac{n \sigma\deg (E)}{2 \Vol(X)} -
n!)(\frac{\omega_{\sigma}^{n+1}}{(n+1)!})\ .\tag 6.4
\endalign
$$
The first result follows from this, Lemma 4.9, and the fact that $\int_{\Bbb
P^1} \omega_{\Bbb P^1} = 1$. The second result then follows by (4.7c) in Lemma
4.9, applied to $\det E$.\qed
\enddemo

We can define $\Cal O_M(1)$\ to be the dual bundle to $\Cal O_M(-1)$, and set
$$\Cal O_M(\pm k)=\Cal O_M(\pm 1)^{\otimes k}\ .$$

\proclaim{Proposition 6.2 (cf. [Ko])}
$$\Cal K_{M/X}^*\equiv\tilde\Cal O(2)=
\Cal O_M(2)\otimes\pi^*(\det E)\ ,$$
where $\det E$\ is the determinant bundle of $E$, and $\pi:M\-->X$\ is as in
the
previous section.
\endproclaim

Combining Propositions 6.1 and 6.2 we thus get
\proclaim{Proposition 6.3 }
\roster
\item
For $k \in \Bbb Z$,
$$\align
\deg_{\sigma}{\Cal O}_M(k) &= k( n! \Vol (X) -
\frac{\deg_{\sigma}\pi^*(\det E)}{2} )\\
&= k(n! \Vol (X) - \frac{n\sigma\deg (E)}{2} )~,\tag 6.5
\endalign$$

\item
$$
\deg_{\sigma} ({\Cal K}_{M/X}^* ) = 2n! \Vol (X)~.\tag 6.6
$$

\endroster
\endproclaim

\noindent\bf Remark 6.4.\rm\  Notice that as a result of (6.6) and (5.11) we
see
that when $M=\Bbb P(E)$\ the parameters $\tau_1,\tau_2$, and $\sigma$\ (in
Theorem 5.1) are related in precisely
the same way as in the case where $M=X\times\Bbb P^1$, i.e. as in  Proposition
3.2.

\heading
\S 7.  Discussion of Results
\endheading

It is clear from the results described above that dimensional
reduction should be thought of in more general terms than
was previously expected. In particular, the role of a global
symmetry in the form of a group action is not as central as
we assumed. On the other hand, what makes the reduction possible
is certainly some kind of geometric order, even if cannot quite
be called a global symmetry. In this section we explore what kind
of geometric structure is evident in our situation, and look at how it
compares to the Garcia-Prada case.

In the case that $M=X\times\Bbb P^1$, there is an action of $SU(2)$\ on $M$\
as
described in [GP]. The action is trivial on $X$\ and via the identification
$\Bbb P^1=SU(2)/U(1)$\ on $\Bbb P^1$. The group orbits are thus of the form
$x\times\Bbb P^1$, and the orbit space is $X$.

In order to exploit this symmetry (i.e. the group action), one must consider
\it invariant\rm\  solutions for whatever equations are to be dimensionally
reduced. In order for such solutions to exist, one needs to be in
an \it equivariant\rm\  setting, and it is for this reason that the bundles
considered in [GP] are $SU(2)$-equivariant bundles.

As soon as the structure of $M$\ as a $\Bbb P^1$-bundle over $X$\ becomes
non-trivial, the $SU(2)$-action is destroyed . Notice however that the group
orbits survive in the form of fibers of the bundle. As we shall see, it is
perhaps more pertinent to describe the fiber $\Bbb P^1$'s as the \it leaves of
a foliation\rm . Indeed, given a foliation on $M$, there is a natural
replacement for the notion of an equivariant bundle, namely that of a \it
foliated bundle\rm\ (cf. [KT]).
Roughly speaking, a foliated bundle has a flat structure along each leaf of
the
foliation. The relation between such bundles and the equivariant bundles
(which
can exist when the foliation comes from a global group action) is made even
clearer if we introduce the holonomy groupoid, $\Cal G_{\Cal F}$, associated
to
the foliation $\Cal F$. The foliated bundles we consider
 can be described as $\Cal
G_{\Cal F}$-equivariant bundles, where an $\Cal G_{\Cal F}$-equivariant bundle
is one which supports an action of $\Cal G_{\Cal F}$\
on its fibers. The action in question is in fact via holonomy transport along
the leaves of the foliation.

This language of foliations provides a convenient framework for
understanding the dimensional reduction we have carried out in \S 5. It is
illuminating even in the case where $M=X\times\Bbb P^1$. In that case the
$SU(2)$-equivariant bundles over $X\times\Bbb P^1$\ can indeed be considered
as
foliated bundles - but in way that is slightly
unexpected.   The novelty is due to the fact that $X\times\Bbb P^1$\ is
foliated in \it two\rm\ ways; in one foliation the leaves are the copies of
$\Bbb P^1$\ in $X\times\Bbb P^1$, and in the other way the leaves are the
copies of $X$.  We will refer to these as the \it fiber foliation\rm\ and the
\it base foliation\rm\ respectively.

\proclaim {Lemma 7.1} The bundles over $X\times\Bbb P^1$\ that are foliated
with respect to the fiber foliation are of the form $\pi_1^*E$, where
$E\-->X$\ is any smooth vector bundle over $X$ and $\pi_1$\ denotes the
projection onto the first factor of $X\times\Bbb P^1$.

The line bundles over $X\times\Bbb P^1$\ that are foliated with respect to the
base foliation are of the form $\pi_2^*H^{\otimes k}$, where  $H\-->\Bbb P^1$\
is the line bundle of degree 1 over $\Bbb P^1$ and $\pi_2$\ denotes the
projection onto the second factor of $X\times\Bbb P^1$.
\endproclaim

Thus the $SU(2)$-equivariant bundles that appear in the dimensional reduction
construction , i.e. the bundles of the form
$V=\pi_1^*E\otimes\pi_2^*H^{\otimes k}$, can be thought of as \it doubly
foliated bundles\rm. The next Proposition shows that it is this doubly
foliated
property that is retained when the trivial fibration $X\times\Bbb P^1$\ is
replaced by a flat projective fibration
$M\-->X$.

\proclaim {Proposition 7.2} Let $M\-->X$\ be a flat projective $PU(2)$-bundle
over $X$, as in (\S 4). Then $M$\ has two natural foliations; one in which the
leaves are the $\Bbb P^1$-fibers of $M\-->X$, and one in which the leaves
are copies of $\tilde X$, the universal cover of $X$. The leaves of
the first foliation intersect the leaves of the second transversally.
\endproclaim

\demo{Proof} This follows immediately form the descriptions (M1-3) of $M$\
given in \S 4.\qed
\enddemo

We will denote the foliation by the $\Bbb P^1$\ fibers by $\Cal F_{\pi}$,
and the other foliation by $\Cal F_{\alpha}$. The next result shows that
the bundles with the structure of extensions as in (4.1) are indeed the
analogues of the bundles considered in [GP].

\proclaim{Proposition 7.3} Let $M\-->X$\ be a flat projective $PU(2)$-bundle
over $X$, as in (\S 4).
\roster
\item There is a one to one correspondence between vector bundles $W$ over $X$
and $\Cal G_{\Cal F_{\pi}}$-equivariant bundles on $M$, given by
$W\-->\pi^*W$.
\item There is a one to one correspondence between $\alpha$-equivariant vector
bundles $V$ over $\Bbb P^1$ and $\Cal G_{\Cal F_{\alpha}}$-equivariant
bundles
on $M$, given by $V\-->\tilde V$, where $\tilde V$\ is the extension as
defined
in \S 4.2
\endroster
\endproclaim

Finally, we remark that $\Bbb P^1$-bundles are almost certainly not the only
objects on which our techniques will work. In particular, there is good reason
to expect that $\Bbb P^1=SU(2)/U(1)$\ can be replaced by any \kler ian
homogeneous space $G/H$. We will return to this question in a later
publication.


\enddocument